\documentclass{amsart}
\usepackage{amssymb}
\usepackage{pdflscape}
\usepackage{amscd}
\usepackage{fancyhdr}

\newtheorem{theorem}{Theorem}[section]

\newtheorem{definition}[theorem]{Definition}

\title{Wild quantum dilogarithm identities}
\author{Markus Reineke}
\begin{document}
\begin{abstract} We exhibit and discuss ``wild'' analogues of the five-term quantum dilogarithm identity. We derive these from the representation theory of quivers, using motivic wall-crossing, the geometricity of motivic Donaldson--Thomas invariants, and special properties of Kronecker moduli.
\end{abstract}
\maketitle
\parindent0pt

\section{Introduction}

The quantum dilogarithm is a $q$-series with many remarkable properties \cite{Z}, including the famous five-term identity \cite{FK}. Cluster algebra theory and wall-crossing of motivic invariants of quivers have led to vast generalizations of such dilogarithm identities \cite{K}.\\[1ex]
In this note, we explore the outer limits of this circle of ideas, by investigating ``wild'' dilogarithm identities, those arising from wild quivers.  For this we use again motivic wall-crossing, interpret the resulting series in terms of motivic Donaldson--Thomas invariants, and use the geometric interpretation of the latter in terms of intersection homology of quiver moduli spaces \cite{MR} to establish very strong positivity properties. In the rank two case, originating from generalized Kronecker quivers, the well-explored and very special symmetries of Kronecker moduli yield many additional explicit properties of the individual terms of our highly infinite dilogarithm identities; see Theorem \ref{mainkronecker}.\\[1ex]
To derive this identity, we collect in Section \ref{s3} the available material on wall-crossing of motivic invariants (see \cite{Moz} for an introduction), and adapt it to the present notation and special setting in Section \ref{s4}. Although similar approaches are used, for example, in the context of the tropical vertex \cite{GP,RW}, it it desirable to state the nature of such wild identities as explicitly as possible.\\[1ex]
{\bf Acknowldegments:} The author would like to thank Daping Weng for several discussions on these wild identities, Bernhard Keller for the opportunity to present them in his seminar, Vladimir Fock as well as Sergey Mozgovoy for valuable comments, and Timm Peerenboom for carefully reading a draft of this text and suggesting several improvements.

\section{Quantum dilogarithm identities}\label{s2}

We first define the quantum dilogarithm, and state the classical five-term identity. The coefficient ring $\mathbb{Q}(q^{1/2})[[x]]$, and in particular the twist by half-powers of $q$, will become natural in the context of motivic Donaldson--Thomas invariants.

\begin{definition} We define the quantum dilogarithm $\Phi(x)\in\mathbb{Q}(q^{1/2})[[x]]$ as 
$$\Phi(x)=\sum_{n\geq 0}\frac{q^{n/2}x^n}{(1-q)\cdot\ldots\cdot(1-q^n)}=$$
 $$=\exp\left(\sum_{n\geq 1}\frac{x^n}{n\cdot (q^{-n/2}-q^{n/2})}\right)=\prod_{n\geq 0}\frac{1}{1-q^{n+1/2}\cdot x}.$$\end{definition}

{\it Remark:} The translation from the present definition to the ones in the literature is straightforward. For example, \cite{FK,V} use the definition $(x;q)_\infty=\Phi(q^{-1/2}x)^{-1}$, and \cite{K} uses $\mathbb{E}(x)=\Phi(-qx)$. The classical Euler dilogarithm \cite{Z} arises as the following limit: $$(q^{-1/2}-q^{1/2})\log\Phi(x)=\sum_{n\geq 1}\frac{x^n\cdot(q^{-1/2}-q^{1/2})}{n\cdot(q^{-n/2}-q^{n/2})}\stackrel{q\rightarrow 1}{\longrightarrow}\sum_{n\geq 1}\frac{x^n}{n^2}= {\rm Li}_2(x)$$

\begin{definition} For a positive integer $m$, we define $\mathbb{Q}(q^{1/2})_{q^m}[[x,y]]$ as the skew formal power series ring with skew commutativity relation $xy=q^myx$.
\end{definition}

We can now formulate the five-term quantum dilogarithm identity:

\begin{theorem}[Schützenberger, Fadeev, Kashaev, Volkov] In $\mathbb{Q}(q^{1/2})_q[[x,y]]$, we have
$$\Phi(x)\Phi(y)=\Phi(y)\Phi(-q^{-1/2}xy)\Phi(x)$$
\end{theorem}

For $m=2$, we will obtain the following identity as a special case of Theorem \ref{mainkronecker}.

\begin{theorem} In $\mathbb{Q}(q^{1/2})_{q^2}[[x,y]]$, we have $$\Phi(x)\Phi(y)=\Phi(y)\Phi(q^{-2}xy^2)\Phi(q^{-6}x^2y^3)\Phi(q^{-12}x^3y^4)\ldots\cdot$$
$$\cdot \Phi(q^{-1/2}xy)^{-1}\Phi(q^{-3/2}xy)^{-1}\cdot$$
$$\cdot\ldots\Phi(q^{-12}x^4y^3)\Phi(q^{-6}x^3y^2)\Phi(q^{-2}x^2y)\Phi(x).$$
\end{theorem}



To make such identities more readable, we will now introduce a shorthand notation, which again will be motivated later by the framework of Donaldson--Thomas invariants:

\begin{definition} For $a,b\geq 0$, define
$$\Phi_{(a,b)}=\Phi({(-1)}^{mab}{q}^{(a^2+b^2-2mab-1)/2}x^ay^b)^{{(-1)}^{a^2+b^2-mab-1}}.$$ 
More generally, for a Laurent polynomial $P(q)=\sum_{k}c_k(-q^{1/2})^k\in\mathbb{Q}[q^{\pm 1/2}]$, define
$$\Phi_{(a,b)}^{{\circ P}}={\prod_{k}}\Phi((-1)^{mab}q^{(a^2+b^2-2mab-1{+k})/2}x^ay^b)^{(-1)^{a^2+b^2-mab-1{+k}}\cdot{c_k}}.$$ 
\end{definition}

Then the previous identities simplify to 
$$\Phi_{(1,0)}\Phi_{(0,1)}=\Phi_{(0,1)}\Phi_{(1,1)}\Phi_{(1,0)}\mbox{ for }
m=1,$$
$$ \Phi_{(1,0)}\Phi_{(0,1)}=\Phi_{(0,1)}\Phi_{(1,2)}\Phi_{(2,3)}\ldots\cdot\Phi_{(1,1)}^{\circ(q+1)}\cdot\ldots\Phi_{(3,2)}\Phi_{(2,1)}\Phi_{(1,0)}$$
for $m=2$.\\[1ex]
In the case $m\geq 3$, we can no longer give explicit identities, but we will obtain a rather complete qualitative description.


%

To formulate this main result, we need two more definitions. We denote by $\sigma$ the operator on $\mathbb{Z}^2$ given by  $\sigma(a,b)=(b,mb-a)$, which generates an infinite dihedral group together with the involution $(a,b)\mapsto (b,a)$. We also define $\mu_\pm={(m\pm\sqrt{m^2-4})}/{2}$, the two roots of the quadratic equation $x^2-mx+1=0$.

\begin{theorem}\label{mainkronecker} In $\mathbb{Q}(q^{1/2})_{q^m}[[x,y]]$ for $m\geq 3$, we have
$$\Phi_{(1,0)}\Phi_{(0,1)}=\Phi_{(0,1)}\Phi_{\sigma(0,1)}\Phi_{\sigma^2(0,1)}\Phi_{\sigma^3(0,1)}\cdot\ldots\cdot$$
$$ \cdot\prod^{\rightarrow}_{{\mu_-\leq a/b\leq \mu_+}\atop{\mbox{\tiny increasing}}}\Phi_{(a,b)}^{\circ P_{(a,b)}}\cdot$$
$$ \cdot\ldots\Phi_{\sigma^{-3}(1,0)}\Phi_{\sigma^{-2}(1,0)}\Phi_{\sigma^{-1}(1,0)}\Phi_{(1,0)},$$ 
where the $P_{(a,b)}$ satisfy the following properties:
\begin{enumerate}
\item {\bf ``Wildness''/Completeness:} We have the non-vanishing property $$P_{(a,b)}\not=0\mbox{ for all }\mu_-\leq a/b\leq \mu_+.$$
\item {\bf Dihedral symmetry:} We have the symmetries
$$P_{(b,a)}=P_{(a,b)},\;\;\;P_{\sigma(a,b)}=P_{(a,b)}.$$
\item {\bf Positivity:} We have
$$P_{(a,b)}\in\mathbb{N}[q],\mbox{ of degree } d=mab-a^2-b^2+1>0.$$
\item {\bf Unimodality:} The polynomial $P_{(a,b)}=\sum_kc_kq^k$ is palindromic and unimodal:
$$c_{d-k}=c_k\mbox{ and }1=c_0\leq c_1\leq \ldots \geq c_{d-1}\geq c_{d}=1.$$
\item {\bf Lowest order terms:}~  $$P_{(1,k)}=\left[{m\atop k}\right]_q.$$
\item {\bf Special value:} ~
$$P_{(k,k)}(1)=\frac{1}{(m-2)k^2}\sum_{d|k}\mu(\frac{k}{d})(-1)^{md+1}{{(m-1)^2d-1}\choose{d}}.$$
\end{enumerate}
\end{theorem}
As a consequence of the dihedral symmetry property, we see that all $P_{(a,b)}$ are determined by those for $a\leq b\leq\frac{m}{2}a$. All properties will follow from interpreting the $P_{(a,b)}$ as the Poincaré polynomials in intersection homology of Kronecker moduli \cite{D}, a class of projective varieties parametrizing certain tuples of matrices up to base change.

\section{Quiver setup}\label{s3}
In this section, we recall the necessary terminology on quiver representations and their moduli spaces, and we formulate the main ingredients for general quiver dilogarithm identities, namely the motivic wall-crossing formula, and the definition and geometricity of motivic Donaldson--Thomas invariants. The reader is referred to \cite{Moz} for a detailed introduction into motivic wall-crossing for quivers, and to \cite{RSmall} for a short summary.\\[1ex]
Let $Q$ be a finite acyclic quiver. We order the set of vertices $Q_0=\{i_1,\ldots,i_n\}$ such that $i_k\rightarrow i_l$ implies $k>l$. The Euler form of $Q$ is the (in general non-symmetric) bilinear form on $\mathbb{Z}Q_0$ given by $\langle\mathbf{d},\mathbf{e}\rangle=\sum_{i\in Q_0}d_ie_i-\sum_{\alpha:i\rightarrow j}d_ie_j$ for $\mathbf{d},\mathbf{e}\in\mathbb{Z}Q_0$.\\[1ex] 
We define the formal quantum affine space $\mathbb{Q}(q^{1/2})_q[[Q_0]]$ as the skew formal power series ring with topological basis $t^\mathbf{d}$ for $\mathbf{d}\in\mathbb{N}Q_0$ and multiplication twisted by the antisymmetrized Euler form
$$t^{\bf d}\cdot t^{\bf e}=(-q^{1/2})^{\langle{\bf d},{\bf e}\rangle-\langle{\bf e},{\bf d}\rangle}t^{{\bf d}+{\bf e}}.$$

We fix linear functions $\Theta,\kappa\in(\mathbb{Z}Q_0)^*$ on $Q$ such that $\kappa({\bf d})>0$ for ${\bf d}\in\mathbb{N}Q_0\setminus 0$, and consider the associated slope function $\mu({\bf d})=\frac{\Theta({\bf d})}{\kappa({\bf d})}$ for $\mathbf{d}\in\mathbb{N}Q_0\setminus 0$. We denote by $\Lambda_a$ the set of all ${\bf d}\in\mathbb{N}Q_0\setminus 0$ of slope $a\in\mathbb{Q}$.\\[1ex]
We define the Grothendieck ring of varieties $K_0({\rm Var}_\mathbb{C})$ as the free abelian group in isomorphism classes of complex algebraic varieties $X$ modulo the ``cut-and-paste'' relation
${ }[X]=[A]+[U]\mbox{ if } A\subset X\mbox{ closed, } U=X\setminus A,$
with product given by $[X]\cdot[Y]=[X\times Y]$. We abbreviate the Lefschetz motive $[\mathbb{A}^1]=:q$.\\[1ex]
We consider the localization $R=K_0({\rm Var}_\mathbb{C})[q^{\pm 1/2}, (1-q^n)^{-1}: n\geq 1]$, and define the formal motivic affine space $R_q[[Q_0]]$ as above, with $R$ replacing the coefficient ring $\mathbb{Q}(q^{1/2})$. In fact, all our computations will happen in the smaller coefficient ring of motives which are rational functions in $q$. Note that the existence of motivic measures such as the virtual Hodge polynomial shows that this subring is isomorphic to a subring of $\mathbb{Q}(q^{1/2})$.\\[1ex]
Given a dimension vector ${\bf d}\in\mathbb{N}Q_0$, we fix $\mathbb{C}$-vector spaces $V_i$ of dimension $d_i$ for $i\in Q_0$. We consider the base change action 
$$\prod_{i\in Q_0}{\rm GL}(V_i)=G_{\bf d}\curvearrowright R_{\bf d}(Q)=\bigoplus_{\alpha:i\rightarrow j}{\rm Hom}(V_i,V_j)$$
given by $$(g_i)_i\cdot(f_\alpha)_\alpha=(g_jf_\alpha g_i^{-1})_{\alpha:i\rightarrow j},$$
whose orbits, by definition, correspond bijectively to the isomorphism classes of complex representations of $Q$ of dimension vector $\mathbf{d}$.\\[1ex]
We denote by $R_{\bf d}^{\mu-{\rm sst}}(Q)\subset R_{\bf d}(Q)$ the open subset of $\mu$-semistable points. Using this notation, we can formulate the motivic wall-crossing formula, which is formally equivalent to the existence of the Harder--Narasimhan filtration \cite{RHN}:

\begin{theorem}[Motivic wall-crossing formula]\label{wcf} In $R_q[[Q_0]]$, we have the identity
$$\Phi(t_1)\cdot\ldots\cdot \Phi(t_n)= \sum_{{\bf d}}(-q^{1/2})^{\langle{\bf d},{\bf d}\rangle}\frac{[R_{\bf d}(Q)]}{[G_{\bf d}]}t^{\bf d}=$$
$$ =\prod^\rightarrow_{a\in\mathbb{Q}\mbox{ \tiny decreasing}}\left(1+\sum_{{\bf d}\in\Lambda_a}(-q^{1/2})^{\langle{\bf d},{\bf d}\rangle}\frac{[R_{\bf d}^{\mu-{\rm sst}}(Q)]}{[G_{\bf d}]}t^{\bf d}\right).$$
\end{theorem} 

Next, we can define motivic Donaldson--Thomas invariants. We generalize the shorthand notation of the previous section to
$$\Phi(x)^{\circ P(q)}=\prod_k\Phi(q^{k/2}x)^{(-1)^kc_k}$$
for $P(q)=\sum_kc_k(-q^{1/2})^k$.

\begin{definition}\cite{KS}\label{dt} Assume that $\langle\_,\_\rangle$ is symmetric on $\Lambda_a$.   Define ${\rm DT}^\mu_{\bf d}(q)\in\mathbb{Q}[q^{\pm 1/2}]$   by factorization in $R_q[[Q_0]]$:
$$1+\sum_{{\bf d}\in\Lambda_a}(-q^{1/2})^{\langle{\bf d},{\bf d}\rangle}\frac{[R_{\bf d}^{\mu-{\rm sst}}(Q)]}{[G_{\bf d}]}t^{\bf d}=\prod_{{\bf d}\in\Lambda_a}\Phi(t^{\bf d})^{\circ {\rm DT}^\mu_{\bf d}(q)}.$$

 The ${\rm DT}^\mu_{\bf d}(q)$ are called the motivic Donaldson--Thomas invariants of the quiver $Q$ with stability $\mu$.
\end{definition}

We remark that the more common definition in terms of the plethystic exponential ${\rm Exp}$ is equivalent to this one since, by definition, $\Phi(x)={\rm Exp}(\frac{x}{q^{-1/2}-q^{1/2}})$.\\[1ex]
The motivic Donaldson--Thomas invariants admit a geometric interpretation in terms of intersection homology of moduli spaces. Namely, we consider the GIT quotient
$$R_{\bf d}^{\mu-{\rm sst}}(Q)//G_{\bf d}=M_{\bf d}^{\mu-{\rm sst}}(Q),$$
 the moduli space of $\mu$-semistable representations of $Q$ of dimension vector ${\bf d}$.\\[1ex] 
It is an irreducible projective normal (typically singular) complex algebraic variety. If ${\bf d}$ is $\mu$-stable, that is, if there exists a $\mu$-stable representation of dimension vector ${\bf d}$, then $\dim M_{\bf d}^{\mu-{\rm sst}}(Q)=1-\langle{\bf d},{\bf d}\rangle$.  In terms of this moduli space, we have the following geometric interpretation of the motivic Donaldson--Thomas invariants \cite{MR}:
 
\begin{theorem}[Geometricity of DT invariants]\label{geom} We have $${\rm DT}_{\bf d}^\mu(q)=(-q^{1/2})^{\langle{\bf d},{\bf d}\rangle-1}\sum_{k\geq 0}\dim{{\rm IH}}^k(M_{\bf d}^{\mu-{\rm sst}}(Q),\mathbb{Q})(-q^{1/2})^k$$
if ${\bf d}$ is $\mu$-stable,   and ${\rm DT}_{\bf d}^\mu(q)=0$ otherwise.
\end{theorem}

\section{Quantum dilogarithm identity for quivers}\label{s4}

To combine the methods prepared in the previous section, we consider quivers and stabilities such that the restriction of $\langle\_,\_\rangle$ is symmetric on all $\Lambda_a$. In particular, this holds if the antisymmetrized Euler form of $Q$ is determined by the stability function, in the sense that
$$\langle{\bf d},{\bf e}\rangle-\langle{\bf e},{\bf d}\rangle=\kappa({\bf d})\Theta({\bf e})-\kappa({\bf e})\Theta({\bf d})\mbox{ for all }{\bf d},{\bf e}\in\mathbb{Z}Q_0,$$ compare \cite[Proposition 5.2]{RSmall}. Besides generalized Kronecker quiver, this property holds, for example, for complete bipartite quivers.

\begin{theorem}\label{quiverdilog} Assume that $\langle\_,\_\rangle$ is symmetric on all $\Lambda_a$. Then we have a factorization
$$\Phi(t_1)\cdot\ldots\cdot\Phi(t_n)=\prod^\rightarrow_{\mu({\bf d})\mbox{ \tiny  decreasing}}\Phi(t^{\bf d})^{\circ(-q^{1/2})^{\langle{\bf d},{\bf d}\rangle-1}\cdot P_{\bf d}(q)}$$
in $R_q[[Q_0]]$, for polynomials $P_{\bf d}(q)$ with the following properties:
\begin{enumerate}
\item {\bf Non-vanishing:} We have $P_{\bf d}(q)\not=0$ if and only if ${\bf d}$ is $\mu$-stable, 
\item {\bf Positivity:} We have $P_{\bf d}\in\mathbb{N}[q]$, of degree $\deg P_{\bf d}(q)=1-\langle{\bf d},{\bf d}\rangle$,
\item {\bf Unimodality:} $P_{\bf d}(q)$ is palindromic and unimodal, 
\item {\bf Simplicity:} We have $P_{\bf d}(q)=1$ if, additionally, $\langle{\bf d},{\bf d}\rangle=1$.

\end{enumerate}
\end{theorem}

This is now readily proved: combining Theorem \ref{wcf}, Definition \ref{dt} and Theorem \ref{geom}, we see that $P_{\bf d}(q)$ is precisely the Poincar\'e polynomial in intersection homology of $M_{\bf d}^{\mu-{\rm sst}}(Q)$ in case ${\bf d}$ is $\mu$-stable, proving the non-vanishing property. Positivity and unimodality are proven in \cite[Corollary 1.2]{MR}. The degree statement follows from the dimension formula for the moduli space, and consequently the simplicity statement follows.\\[1ex]
To derive Theorem \ref{mainkronecker}, we consider the $m$-arrow Kronecker quiver $K_m=1\stackrel{(m)}{\leftarrow}2$, with stability function given by $\Theta({\bf d})=m\cdot(d_2-d_1)$ and $\kappa({\bf d})=d_1+d_2$.\\[1ex] 
We identify the variables $t_1=x$ and $t_2=y$, leading to an identification of $\mathbb{Q}(q^{1/2})_q[[(K_m)_0]]$ with $\mathbb{Q}(q^{1/2})_{q^m}[[x,y]]$, such that $t^{\bf d}=(-q^{1/2})^{md_1d_2}x^{d_1}y^{d_2}$ and $$\Phi(t^{{\bf d}})^{\circ (-q^{1/2})^{\langle{\bf d},{\bf d}\rangle-1}\cdot P(q)}=\Phi_{(a,b)}^{\circ P(q)}$$
for ${\bf d}=(a,b)$.\\[1ex]
Theorem \Ref{quiverdilog} then provides a factorization in $\mathbb{Q}(q^{1/2})_{q^m}[[x,y]]$ of the form
$$\Phi_{(1,0)}\Phi_{(0,1)}=\prod^\rightarrow_{a/b\mbox{ \tiny increasing}}\Phi_{(a,b)}^{\circ P_{(a,b)}},$$ 
such that
$$P_{(a,b)}(q)=\sum_{k\geq 0}\dim{\rm IH}^{2k}(K_{(a,b)},\mathbb{Q})q^k,$$
where
$$  K_{(a,b)}=M_{a\times b}(\mathbb{C})^m_{\rm sst}//{\rm GL}_a(\mathbb{C})\times{\rm GL}_b(\mathbb{C})$$
are the Kronecker moduli of \cite{D}.\\[1ex] 
All remaining properties in Theorem \ref{mainkronecker} now follow from properties of these Kronecker moduli spaces. Using the analysis of \cite[Section 4]{Scho}, we see that a dimension vector $(a,b)$ is $\mu$-stable if and only if $\langle(a,b),(a,b)\rangle\leq 1$. In case $\langle(a,b),(a,b)\rangle=1$, the dimension vector is a real root, and the moduli space $K_{(a,b)}$ reduces to a point, thus $P_{(a,b)}=1$. Otherwise, we have $\langle(a,b),(a,b)\rangle\leq 0$, which translates to $\mu_-\leq a/b\leq \mu_+$. This proves the completeness statement of the theorem. By linear duality, we have $K_{(a,b)}\simeq K_{(b,a)}$. Moreover, reflection functors induce isomorphisms
$K_{(a,b)}\simeq K_{(b,mb-a)}$ by \cite[Proposition 4.3]{weist}. This establishes the dihedral symmetry. 
The Kronecker moduli $K_{(1,k)}$ identify with the Grassmannians ${\rm Gr}_k^m$, determining the lowest order terms in the factorization. Finally, the special value $P_{(k,k)}(1)$ is computed in \cite[Theorem 5.2]{RFun}.





\end{document}